\documentclass{article}
\usepackage{amsmath,amssymb,amsthm}
\usepackage{graphics}

\newtheorem{theorem}{Theorem}[section]
\newtheorem{lemma}{Lemma}[section]

\newcommand{\dist}{\mathop{\rm dist}\nolimits}
\newcommand{\tr}{\mathop{\rm Tr}\nolimits}
\newcommand{\rk}{\mathop{\rm rank}\nolimits}
\newcommand{\psd}{\mathop{\rm PD}\nolimits}

\newtheorem{defn}{Definition}[section]

\title {Positive definite functions in distance geometry}

\author {Oleg R. Musin \thanks{The author partially supported by NSF and NSA grants.}}

\begin{document}
\date{}
\maketitle

\begin{abstract} I. J. Schoenberg  proved that a function  is positive definite  in the unit sphere
 if and only if this function is a nonnegative  linear combination of Gegenbauer  polynomials.
 This fact play a crucial role in Delsarte's method for finding bounds
 for the density of sphere packings on spheres and Euclidean spaces.

One of the most excited applications of Delsarte's method is  a solution of the kissing number problem in dimensions 8 and 24. However, 8 and 24 are the only dimensions in which this method gives a precise result. For other dimensions (for instance, three and four) the upper bounds exceed the lower. We have found an extension of the Delsarte method
that allows to solve the kissing number problem (as well as the one-sided kissing number problem) in dimensions three and four.

In this paper we also will discuss the maximal cardinalities of spherical two-distance sets.  Using the so-called  polynomial method and Delsarte's method these cardinalities can be determined for all dimensions $n<40$.

Recently, were found extensions of Schoenberg's theorem for multivariate positive-definite functions. Using these extensions  and semidefinite programming can be improved some upper bounds for spherical codes.
\end{abstract}

\section{Introduction}

Let $M$ be a metric space with a distance function $\tau.$ A real continuous function $f(t)$  
is said to be positive definite (p.d.) in $M$ if
for arbitrary points $p_1,\ldots,p_r$ in $M$, real variables $x_1,\ldots,x_r$, and arbitrary $r$  we have
$$
\sum\limits_{i,j=1}^r {f(t_{ij})\,x_ix_j}\ge 0, \quad t_{ij}=\tau(p_i,p_j),
$$
or equivalently,  the matrix $\bigl(f(t_{ij})\bigr)\succeq0$, where  the sign $\succeq 0$ stands for: ``is  positive semidefinite".

Let ${\Bbb S}\sp {n-1}$ denote the unit sphere in ${\Bbb R}\sp n$, and let $\varphi_{ij}$ denote the angular distance between points $p_i, p_j.$  Schoenberg \cite{Scho} proved that:\\
\centerline{\em $f(\cos\varphi)$ is p.d. in ${\Bbb S}\sp {n-1}$ if and only if
$\; f(t)=\sum_{k=0}^\infty{f_kG_k^{(n)}(t)}$ with all $f_k\ge 0$.}
\\
 Here $G_k^{(n)}(t)$
are the Gegenbauer   polynomials.

Schoenberg's theorem has been generalized by Bochner \cite{Boc} to more general spaces. Namely, the following fact holds: {\em $f$ is p.d. in a 2-point-homogenous space $M$ if and only if $f(t)$ is a nonnegative linear combination   of the zonal spherical functions  $\Phi_k(t)$} (see details in \cite[Th. 2]{Kab}, \cite[Chapter 9]{CS}).

Note that
the Bochner - Schoenberg theorem is widely used in coding theory and discrete geometry for finding bounds for error-correcting codes, constant weight codes, spherical codes, sphere packings and other packing problems in 2-point-homogeneous spaces (see \cite{CS,  Kab, Mus2,Mus3, Mus4, PZ} and many others).

The paper is organized as follows:

Section 2 recalls definitions of Gegenbauer polynomials and considers Delsarte's method for spherical codes.

Section 3 discusses applications of Delsarte's method for the kissing number problem.
One of the most excited applications of Delsarte's method is  a solution of the kissing number problem in dimensions 8 and 24. However, 8 and 24 are the only dimensions in which this method gives a precise result. For other dimensions (for instance, three and four) the upper bounds exceed the lower. We have found an extension of the Delsarte method that allows to solve the kissing number problem (as well as the one-sided kissing number problem) in dimensions three and four.

Section 4 discusses maximal cardinalities  of spherical  two-distance sets. Using the so-called  polynomial method and Delsarte's method these cardinalities can be determined for all dimensions $n<40$.

Section 5 considers Sylvester's theorem and semidefinite programming (SDP) bounds for codes. Delsarte's method and its extensions allow to consider the upper bound problem for codes in 2-point-homogeneous spaces as a
linear programming problem with perhaps infinitely many variables, which are the distance distribution. We show that using as variables
power sums of distances this problem can be considered as a finite semidefinite programming  problem. This method allows to improve some
linear programming upper bounds. 

Section 6 discusses an application of the extended Schoenberg's theorem  for multivariate
 Gegenbauer polynomials. This extension derives new positive semidefinite constraints for the distance distribution which can be applied for spherical codes.

\section{Delsarte's method}

\noindent{\bf 2-A. The Gegenbauer polynomials.}
Let us recall definitions of Gegenbauer polynomials. Let polynomials $C_k^{(n)}(t)$ are defined by the expansion
$$(1-2rt+r^2)^{(2-n)/2} = \sum\limits_{k=0}\limits^{\infty}r^kC_k^{(n)}(t).$$
Then the polynomials $G_k^{(n)}(t): = C_k^{(n)}(t)/C_k^{(n)}(1)$ are called {\it Gegenbauer} or {\it ultraspherical} polynomials. (So the normalization of $G_k^{(n)}$ is determined by the condition $G_k^{(n)}(1)=1.$)
Also the Gegenbauer polynomials $G_k^{(n)}$ can be defined by the recurrence formula:
$$G_0^{(n)}=1,\;\; G_1^{(n)}=t,\; \ldots,\; G_k^{(n)}=\frac {(2k+n-4)\,t\,G_{k-1}^{(n)}-(k-1)\,G_{k-2}^{(n)}} {k+n-3}$$

Note that for any even $k\ge0$ (resp. odd)  $G_k^{(n)}(t)$ is even (resp. odd). Therefore,
$G_{2k}$ and $G_{2\ell+1}$ are orthogonal on $[-1,1]$.
Moreover, all  polynomials $G_k^{(n)}$ are orthogonal on $[-1,1]$ with respect to the weight
 function $(1-t^2)^{(n-3)/2}$:
$$
\int_{-1}^1G_k^{(n)}(t)\,G_\ell^{(n)}(t)\,(1-t^2)^{(n-3)/2}\,dt=0, \quad k\ne\ell. $$

Recall the { addition theorem for Gegenbauer polynomials}:
$$
G_k^{(n)}(\cos{\theta_1}\cos{\theta_2}+\sin{\theta_1}\sin{\theta_2}\cos{\varphi})$$
$$=\sum\limits_{s=0}^k c_{nks}\,G_{k-s}^{(n+2s)}(\cos{\theta_1})\,G_{k-s}^{(n+2s)}(\cos{\theta_2})\,(\sin{\theta_1})^s\,(\sin{\theta_2})^s\,G_s^{(n-1)}(\cos{\varphi}),
$$
where $c_{nks}$ are  {\em positive} coefficients whose values of  no concern here  (see \cite{Car,Erd}).

\medskip

\noindent{\bf 2-B. Schoenberg's theorem.} Using the addition theorem for Gegenbauer polynomials Schoenberg \cite{Scho} proved the following theorem:
\begin{theorem} A continuous real function  $f(\cos\varphi)$ is positive definite in ${\Bbb S}\sp {n-1}$ if and only if $f$ is a nonnegative linear combination of the Gegenbauer polynomials, i.e.
$\; f(t)=\sum_{k=0}^\infty{f_kG_k^{(n)}(t)}$ with all $f_k\ge 0$.
\end{theorem}

This theorem also can be proved using the addition theorem for harmonic polynomials (see details in \cite{PZ}).

\medskip

\noindent{\bf 2-C. Delsarte's inequality.}
Let $\{p_1, p_2,\ldots, p_M\}$ be any finite subset of the unit sphere ${\Bbb S}^{n-1}.$ By $\varphi_{ij}=\dist(p_i,p_j)$ we denote the spherical (angular) distance between  $p_i,\,  p_j.$ Clearly,  $ \cos{\varphi_{ij}}=\langle p_i,p_j\rangle.$

If a symmetric matrix is positive semidefinite,  then the sum of all its entries is nonnegative. Schoenberg's theorem implies that
 the matrix $\big(G_k^{(n)}(t_{ij})\big)$ is positive semidefinite, where $t_{ij}:=\cos{\varphi_{ij}}, \; $ Then
$$\sum\limits_{i=1}^M\sum\limits_{j=1}^M {{G_k^{(n)}(t_{ij})}} \ge 0 \eqno (2.1)$$

Suppose a continuous functions $f:[-1,1]\to {\bf R}$ is p.d. in ${\Bbb S}\sp {n-1}$. Then
$$
f(t)=\sum\limits_{k=0}^\infty {f_kG_k^{(n)}(t)}
$$
with all $f_k\ge 0$.
Let  $$S(X)=S_f(X):=\sum\limits_{i=1}^M\sum\limits_{j=1}^M{f(t_{ij})}.$$ Using $(2.1),$ we get
$$S(X)=
\sum\limits_{k=0}^\infty f_k\left(\sum\limits_{i=1}^M\sum\limits_{j=1}^M {G_k^{(n)}(t_{ij})}\right)\ge
\sum\limits_{i=1}^M\sum\limits_{j=1}^M {f_0G_0^{(n)}(t_{ij})} =   f_0M^2.$$
Thus
$$S_f(X)\ge f_0M^2. \eqno (2.2)$$

\medskip

\noindent{\bf 2-D. Delsarte's bound.} We say that
$X=\{p_1,\ldots,p_M\}\subset {\Bbb S}^{n-1}$
is a {\it spherical $\psi$-code,} where $0<\psi<\pi$, if for all $i\neq j,$ $t_{ij}=\cos{\phi_{ij}} \le z:=\cos{\psi},$ i.e. $t_{ij}\in [-1,z]$. In other words, the angular separation between distinct points from $X$ is at least $\psi$.
{Denote by  $A(n,\psi)$  the maximal size of a $\psi$-code in ${\Bbb S}^{n-1}$.}

\begin{theorem}[\cite{Del1,Del2,Kab}] Let a continuous function $f:[-1,1]\to {\bf R}$ be p.d. in ${\Bbb S}\sp {n-1}$. Let
$\; f(t)\le 0\; $  for all $t\in [-1,\cos{\psi}].\; $
Then
$$A(n,\psi) \le \frac {f(1)}{f_0}.$$
\end{theorem}
\begin{proof} Let $X=\{p_1,\ldots,p_M\}\subset {\Bbb S}^{n-1}$ be a spherical $\psi$-code. Clearly, $f(t_{ii})=f(1)$.
By assumptions we have $\; f(t_{ij})\le 0\; $ for all $\; i\neq j$. Therefore
 $$S_f(X)=Mf(1)+2f(t_{12})+\ldots+2f(t_{M-1,M}) \le Mf(1).$$  If we combine this with $(2.2),$ then we get
 $M \le f(1)/f_0.$
 \end{proof}

\section{The kissing problem}

\noindent{\bf 3-A. The kissing number problem.}
The kissing number problem asks for the maximal number $k(n)$ of equal size nonoverlapping
spheres  in $n$-dimensional space  that can touch another sphere of the same size. In other words, $k(n)=A(n,\pi/3)$, i.e. $k(n)$ is the maximal size
of a spherical $\pi/3$-code of length (dimension) $n$.

This problem in dimension three was the subject of a famous discussion between Isaac Newton
and David Gregory in 1694. In three dimensions the problem was finally solved only in 1953 by Sch\"utte and van der Waerden \cite{SvdW2}.

In 1979
Levenshtein \cite{Lev2}, and independently Odlyzko and  Sloane \cite{OdS} (= \cite[Chap.13]{CS}), using Delsarte's method, have proved that $k(8)=240$, and $k(24)=196560$.  Moreover, Bannai and Sloane \cite{BS} (= \cite[Chapter 14]{CS})  proved that the maximal kissing  arrangements in these dimensions are unique up to isometry.
However, $n=8, 24$ are the only dimensions in which this method gives a precise result. For other dimensions (for instance, $n=3, 4$) the upper bounds exceed the lower.

\medskip

\noindent{\bf 3-B. The kissing problem in dimensions 8 and 24.}
The proofs in \cite{Lev2,OdS} that $k(8)=240$ and $k(24)=196560$ are surprisingly short, clean, and technically easier than all known proofs in three dimensions. Indeed, let
$$
f_8(t)=(t-1/2)\,t^2(t+1/2)^2(t+1)=\sum\limits_{k=0}^6{f_k^{(8)}G_k^{(8)}(t)},
$$
and
$$
f_{24}(t)=(t-1/2)(t-1/4)^2\,t^2(t+1/4)^2(t+1/2)^2(t+1)=\sum\limits_{k=0}^{10}{f_k^{(24)}G_k^{(24)}(t)}.
$$
Since all $f_k^{(8)}\ge 0, \; f_k^{(24)}\ge 0$ and $f_8(t)\le0,\; f_{24}(t)\le 0$ for all $t\in [-1,1/2]$ Theorem 2.2 yields
$$
k(8)\le \frac{f_8(1)}{f_0^{(8)}}=240, \quad k(24)\le \frac{f_{24}(1)}{f_0^{(24)}}=196560.
$$
For $n=8, \, 24$ the minimal vectors in sphere packings $E_8$ and Leech lattices give these kissing numbers. Thus $k(8)=240,$ and $k(24)=196560.$

\medskip

\noindent{\bf 3-C. The kissing problem in four dimensions.}
It is not hard to see that $k(4)\ge 24$.
Indeed, the unit sphere in ${\Bbb R}^4$ centered at $(0,0,0,0)$ has 24 unit spheres around it, centered at the points $(\pm\sqrt{2},\pm\sqrt{2},0,0)$, with any choice of signs and any ordering of the coordinates. The convex hull of these 24 points yields a famous 4-dimensional regular polytope - the ``24-cell". Its facets are 24 regular octahedra.

Let $f_{OS}(t)=f_0+f_1G_1^{(4)}(t)+\ldots+f_9G_9^{(4)}(t),\; $  where
$f_0=1, \; f_1=3.6181,\; f_2=6.1156,\; f_3=7.0393,\; f_4=5.0199,\;
f_5=2.313, \;  f_6=f_7=f_8=0, \; f_9=0.4525.$ {This polynomial was applied by Odlyzko and Sloane \cite{OdS} to prove that $k(4)\le 25.$ Since $f_{OS}(t)\le0$ for $t\in[-1,1/2]$, Delsarte's bound gives $$k(4)=A(4,\pi/3)\le  f_{OS}(1)/f_0=f_{OS}(1)\approx 25.5584.$$ Thus, $\, k(4)\le 25.$}

Note that Arestov and Babenko \cite{AB1} proved that
the bound $\; k(4)\le25\; $ cannot be improved using Delsarte's method.

Let
$$f_4(t): = \frac{1344}{25}\,t^9 - \frac{2688}{25}\,t^7 + \frac{1764}{25}\,t^5 + \frac{2048}{125}\,t^4 - \frac{1229}{125}\,t^3 - \frac{516}{125}\,t^2 -
 \frac{217}{500}\,t - \frac{2}{125}$$
In \cite{Mus2} we proved that  $k(4)=24$. This proof is based on the  following two lemmas:

\begin{lemma} { Let $X=\{x_1,\ldots,x_M\}$ be points in the unit sphere ${\bf S}^3$. Then
$$
S(X)=\sum\limits_{i=1}^M\sum\limits_{j=1}^M {f_4(\langle x_i,x_j\rangle)}\ge M^2.
$$
  }
\end{lemma}
\begin{proof}
The expansion of $f_4$ in terms of $U_k=G_k^{(4)}$ is
$$
f_4 = \sum\limits_{i=0}^9 {f_i^{(4)}U_i} = U_0 + 2\,U_1 + \frac{153}{25}\,U_2 + \frac{871}{250}\,U_3 +
\frac{128}{25}\,U_4
+\frac{21}{20}\,U_9.
$$
We see that all $f_i^{(4)}\ge0$ and $f_0^{(4)}=1$. So  Lemma 3.1 follows from $(2.2)$.
\end{proof}

\medskip

 \begin{lemma}{ Suppose $X=\{x_1,\ldots,x_M\}$ is a subset of ${\bf S}^3$ such that the angular separation between any two distinct points $x_i, x_j$ is at least $\pi/3$. Then
$$
S(X)=\sum\limits_{i=1}^M\sum\limits_{j=1}^M {f_4(\langle x_i,x_j\rangle)} < 25M.
$$
  }
\end{lemma}

It's not easy to prove this lemma. A proof is given in \cite[Sections 4,5,6]{Mus2}.

\begin{theorem}
$\quad k(4)=24.$
\end{theorem}
\begin{proof} Let $X$ be a spherical $\pi/3$-code in ${\bf S}^3$ with $M=k(4)$ points.
Then $X$ satisfies the assumptions in Lemmas 3.1, 3.2.  Therefore, $M^2\le S(X) <25M.$
From this $M<25$ follows, i.e. $M\le 24$. From the other side we have $k(4)\ge 24$,
showing that $M=k(4)=24.$
\end{proof}

\medskip

\noindent{\bf 3-D. The kissing problem in three dimensions.}
Our extension of the Delsarte method can be applied to other dimensions and spherical $\psi$-codes.
 The most interesting application is a new proof for the Newton-Gregory problem, $k(3)<13.$ In dimension three all computations  are technically  much more easier than for $n=4$ (see \cite{Mus13}).

Let
$$f_3(t) = \frac{2431}{80}t^9 - \frac{1287}{20}t^7 + \frac{18333}{400}t^5 + \frac{343}{40}t^4 - \frac{83}{10}t^3 - \frac{213}{100}t^2+\frac{t}{10} - \frac{1}{200}. $$
Then for any $M$-point kissing arrangement $X$ we have $S_{f_3}(X)\le 12.88M$ (see details in \cite{MusJ,Mus13}).
The expansion of $f_3$ in terms of Legendre polynomials $P_k=G_k^{(3)}$ is
$$f_3 = P_0 + 1.6P_1 + 3.48P_2 + 1.65P_3 + 1.96P_4 + 0.1P_5 + 0.32P_9.$$
Since $f_0^{(3)}=1,\; f_i^{(3)}\ge 0,$ we have $S_{f_3}(X)\ge M^2$. Thus, $  k(3)\le 12.88<13.$

\medskip

\noindent{\bf 3-E. The one-sided kissing problem in four dimensions.}
Let $H$ be a closed half-space of ${\bf R}^{n}$. Suppose $S$ is a unit sphere in $H$ that touches the supporting hyperplane of $H$. The one-sided kissing number $B(n)$ is the maximal number of unit nonoverlapping spheres in $H$ that can touch $S$.

If  nonoverlapping  unit spheres kiss (touch) the unit sphere $S$ in $H\subset{\bf R}^n$, then the set of kissing points
is an arrangement on the closed hemisphere $S_+$ of $S$ such that the (Euclidean) distance between any two points is at least 1. So the one-sided kissing number problem can be stated in another way: How many points can be placed on the surface of $S_+$  so that the angular separation between any two points is at least $\pi/3$? In other words, $B(n)$ is the maximal cardinality of a $\pi/3$-code on the hemisphere $S_+$.

Clearly, $B(2)=4$. It is not hard to prove that $B(3)=9.$
 Using  extensions of Delsarte's method we proved that $B(4)=18$ (see \cite{Mus3} for a proof and references). Recently have been  obtained several new upper bounds for the one-sided kissing numbers \cite{BM,Mus4,BV2}.

\section{Spherical two-distance sets}

\noindent{\bf 4-A. Two-distance sets.} A set $S$ in Euclidean space ${\Bbb R}^n$ is called a {\it two-distance set}, if there are two distances $c$ and $d$, and the distances between pairs of points of $S$ are either $c$ or $d$.
If a two-distance set $S$ lies in the unit sphere ${\Bbb S}^{n-1}$, then $S$ is called {\it spherical two-distance set.} In other words, $S$ is a set of unit vectors, there are  two real numbers $a$ and $b$,  $-1\le a,b<1$, and inner products  of  distinct vectors of $S$ are either $a$ or $b$.

The ratios of distances of two-distance sets are quite restrictive. Namely, Larman, Rogers, and Seidel \cite {LRS} have proved the following fact:{ if the cardinality of a two-distance set $S$ in ${\Bbb R}^n,\;$  with distances $c$ and $d, \; c< d,$ is greater than  $2n+3$,  then the ratio $c^2/d^{\,2}$ equals $(k-1)/k$ for an integer   $k$ with
$$
2\le k \le \frac{1+\sqrt{2n}}{2}\,.
$$
}

Einhorn and Schoenberg \cite{ES} proved that there are finitely many two-distance sets $S$ in ${\Bbb R}^n$ with cardinality $|S|\ge n+2$. Delsarte, Goethals, and Seidel \cite{Del2} proved that the largest cardinality of spherical two-distance sets in ${\Bbb R}^n$ (we denote it by $g(n)$) is bounded by $n(n+3)/2$, i.e.,
$$
g(n)\le \frac{n(n+3)}{2}\,.
$$
 Moreover, they give examples of spherical two-distance sets  with $n(n+3)/2$ points for $n=2, 6, 22$. {(Therefore, in these dimensions we have equality $g(n)=n(n+3)/2$.)  }
Blockhuis \cite{Blo1} showed that the cardinality of (Euclidean) two-distance sets in ${\Bbb R}^n$ does not exceed $(n+1)(n+2)/2$.

 The standard unit vectors $e_1,\ldots,e_{n+1}$ form an orthogonal basis of ${\Bbb R}^{n+1}$. Denote by $\Delta_n$ the  regular simplex with vertices $2e_1,\ldots,2e_{n+1}$. Let $\Lambda_n$ be the set of points $e_i+e_j, \; 1\le i<j\le n+1.$ Since $\Lambda_n$ lies in the hyperplane $\sum {x_k}=2$, we see that $\Lambda_n$ represents a spherical two-distance set in ${\Bbb R}^n$. On the other hand, $\Lambda_n$ is the set of mid-points of the edges of $\Delta_n$. Thus,
$$
g(n)\ge |\Lambda_n|=\frac{n(n+1)}{2}\,.
$$

For $n<7$ the largest cardinality of Euclidean two-distance sets is $g(n)$, where $g(2)=5, \; g(3)=6, \; g(4)=10, \; g(5)=16$, and $g(6)=27$ (see \cite{Lis}). However, for $n=7,8\; $ Lison\v{e}k \cite{Lis} discovered non-spherical maximal two-distance sets of the cardinality 29 and 45 respectively.

\medskip

\noindent{\bf 4-B. Spherical two-distance sets with $a+b\ge0$.} In \cite{Mus6} using the polynomial method we proved the following fact:
\begin{theorem} Let $S$ be a spherical two-distance set in ${\Bbb R}^n$ with inner products $a$ and $b$. If $a+b\ge0$, then
$$
|S|\le \frac{n(n+1)}{2}\,.
$$
\end{theorem}

Recently, Nozaki \cite{Noz} extended this theorem for spherical $d$-distance sets.

\begin{theorem}
Let $S$ be a spherical $d$-distance set in ${\Bbb R}^n$ with inner products $a_1,\ldots,a_d$.
 Let
 $$
\prod\limits_{k=1}^d{(t-a_k)}=\sum\limits_{k=0}^d{f_kG_k^{(n)}(t)}.
  $$
  Then
$$
|S|\le \sum\limits_{k: f_k>0}{h_k},
$$
where
$$
h_k={n+k-2\choose k}+{n+k-3\choose k-1}.
$$
\end{theorem}

\medskip

\noindent{\bf 4-C. Delsarte's method for spherical two-distance sets.}
Let $S$ be a spherical two-distance set in ${\Bbb R}^n$ with inner products $a$ and $b$, where $a> b$. Let $c=\sqrt{2-2a}, \; d=\sqrt{2-2b}$. Then $c$ and $d$ are the Euclidean distances of $S$.

Let
$$
b_k(a)=\frac{ka-1}{k-1}\,.
$$

If $k$ is defined by the equation: $b_k(a)=b$, then $(k-1)/k=c^2/d^{\,2}$.
Therefore,
if $|S|>2n+3$,   then $k$ is an integer number and $k\in\{2,\ldots,K(n)\}$ \cite{LRS}. Here, $K(n) = \lfloor \frac{1+\sqrt{2n}}{2} \rfloor$.

Consider the case  $a+b_k(a)<0$. Since $b_k(a)\ge-1$, we have
$$
a\in I_k:=\left[\frac{2-k}{k}\,, \frac{1}{2k-1}\right).
$$

Therefore, for fixed $n,\; k\in\{2,\ldots,K(n)\},$ and $a\in I_k$ we have spherical codes with two inner products $a$ and $b_k(a)$. The maximal cardinality of these codes can be bounded by Delsarte's method  (see details in  \cite{Mus6}). Since   for  $6<n<40, \; n\ne 22,23$, Delsarte's bounds are not greater than $n(n+1)/2$,  we  have that $g(n)=n(n+1)/2$. For $n=23$ we obtain $g(23)\le 277$. But $g(23)\ge 276$. This proves the following theorem:

\begin{theorem} If $\; 6<n<22\; $ or $\; 23<n<40$, then
$$
g(n)=\frac{n(n+1)}{2}\,.
$$
For $n=23$ we have
$$
g(23)=276 \; \mbox{ or } \; 277.
$$
\end{theorem}

The case $n=23$ is very interesting. In this dimension the maximal number of equiangular lines (or equivalently, the maximal cardinality of a two-distance set with $a+b=0$) is $276$ \cite{LeS}. On the other hand, $|\Lambda_{23}|=276$. So in $23$ dimensions we have  two very different two-distance sets with $276$ points.

Note that for $n=23$ the Delsarte's bound is $\approx 277.095$. So this numerical bound is not far from $277$. Perhaps stronger tools, such as semidefinite programming bounds, are needed here to prove that $g(23)=276.$

Our numerical calculations show that the barrier $n=40$ is in fact fundamental: Delsarte's bounds are incapable of resolving the $n\ge40,\; k=2$ case. That means a new idea is required to deal with $n\ge40.$

\section{Sylvester's theorem and SDP bounds for codes}

\noindent {\bf 5-A. 2-point-homogeneous spaces.}
We say that a $G$-space ${\bf M}$ is a {\em 2-point-homogeneous space} if (i) ${\bf M}$ is a metric space with a distance  $\rho$ defined on it;
(ii) ${\bf M}$ is a set on which a group $G$ acts;  (iii) $\rho$ is strongly invariant under $G$, i.e. for $x,x',y,y' \in {\bf M}$ with $\rho(x,y)=\rho(x',y')$  there is an element $g\in G$ such that $g(x)=x'$ and $g(y)=y'$.

These assumptions are quite restrictive. In fact, if $G$ is infinite and ${\bf M}$ is a compact space, then  Wang \cite{Wang} has proved that ${\bf M}$ is a sphere; real, complex or quaternionic projective space; or the Cayley projective plane. However, the finite  2-point-homogeneous spaces have note yet been completely classified (see for the most important examples and references  \cite{CS}).

With any  2-point-homogeneous space ${\bf M}$ and an integer number $k\ge 0$ are associated the {\em zonal spherical function} $\Phi_k(t)$ such that $\{\Phi_k(t)\}_{k=0,1,2,\ldots}$  are orthogonal on
$T:=\{\tau(x,y): x, y\in {\bf M}\},$ where $\tau$ is the certain function in $\rho$ ( i.e
$\tau(x,y)=F(\rho(x,y))$) defined by ${\bf M}$.
For all continuous compact ${\bf M}$  and for all currently known finite cases: {\em $\Phi_k(t)$
 is a polynomial of degree } $k$. The normalization is given by the rule:
$\Phi_k(\tau_0)=1$, where $\tau_0:=\tau(x,x)$. Then  $\Phi_0(t)=1.$

Note that if ${\bf M}$ is a Hamming space ${\bf F}_2^n$ with $\tau(x,y)=\rho(x,y)$=Hamming distance, then $\tau_0=0$.  Here  $\Phi_k(t)$ is the Krawtchouk polynomial $K_k(t,n)$.
  Consider the case ${\bf M}$ = unit sphere ${\Bbb S}^{n-1}\subset {\Bbb R}^n$ with $\tau(x,y)=\cos{\rho(x,y)}=\langle x,y\rangle$, where  $\rho(x,y)$ is the angular distance between $x$ and $y$. Then
the corresponding zonal spherical function $\Phi_k(t)$ is the Gegenbauer  polynomial $G_k^{(n)}(t)$.

\medskip

\noindent {\bf 5-B. The Bochner-Schoenberg theorem.}
The main property for zonal spherical functions is called ``positive-definite degenerate kernels" or ``p.d.k"  \cite{CS}.
This property first was discovered by Bochner \cite{Boc} (general spaces)  and independently for spherical spaces by Schoenberg \cite{Scho}.

Now we explain what the p.d.k. property means for finite subsets in ${\bf M}$.
\begin{theorem}[\cite{Boc,Scho,Kab}] Let ${\bf M}$ be a 2-point-homogeneous space.
Then for any integer $k\ge 0$ and for any finite  $C=\{x_i\}\subset {\bf M}$ the matrix $\left(\Phi_k(\tau(x_i,x_j))\right)$ is positive semidefinite.
\end{theorem}

This theorem implies the fact that plays a crucial role for the linear programming  bounds.
For any positive semidefinite matrix the sum of its entries is nonnegative. Then
\begin{theorem}[\cite{Del1,Del2,Kab,OdS}] For any finite  $C=\{x_i\}\subset {\bf M}$ we have
$$
\sum\limits_{i=1}^{|C|}\sum\limits_{j=1}^{|C|} {\Phi_k(\tau(x_i,x_j))}\ge 0.
$$
\end{theorem}

Since 
$\Phi_k(\tau(x_i,x_i))=\Phi_k(\tau_0)=1$, this theorem implies
$${\frac{1}{|C|} \sum\limits_{i,j:i\ne j} {\Phi_k(\tau(x_i,x_j))}\ge -1}, \; \, k=0, 1, 2, \ldots \, . \eqno (5.1)$$

\medskip

\noindent {\bf 5-C. The linear programming bounds.}
Let $S$ be a fixed subset of $T$. We say that a finite subset $C\subset {\bf M}$ is an {\em $S$-code} if $\tau(x,y)\in S$ for all $x, y \in C, x\ne y.$ The largest cardinality $|C|$ of an $S$-code will be denoted by $A({\bf M},S)$.


 The {\em distance distribution} $\{\alpha_t\}$ of $C$ is defined by
$$
\alpha_t:=\frac{1}{|C|}(\mbox{number of ordered pairs }\, x, y\in C \, \mbox{ with }\, \tau(x,y)=t).
$$
We obviously have
$$
\alpha_{\tau_0}=1, \quad \alpha_t\ge 0, \; t\in T, \quad \sum\limits_{t\in T}{\alpha_t}=|C|. \eqno (5.2)
$$

$(5.1)$ and $(5.2)$ make it possible to regard the problem of bounding $A({\bf M},S)$ as a linear programming  problem:

\medskip

\noindent{\em Primal problem (LPP):} Choose a natural number $s$, a subset $\{\tau_1, \ldots, \tau_s\}$ of
  $S$, and real numbers $\alpha_{\tau_1}, \ldots, \alpha_{\tau_s}$ so as to
$$
\mbox{maximize }\; \alpha_{\tau_1} + \ldots + \alpha_{\tau_s}
$$
subject to
$$
\alpha_{\tau_i}\ge 0, \; i=1, \ldots, s, \quad \sum\limits_{i=1}^s {\alpha_{\tau_i}\Phi_k(\tau_i)}\ge -1, \; k=0,1, \ldots .
$$

This is a linear programming problem with perhaps infinitely many unknowns $\alpha_{t}$ and constraints $(1), (2)$.
If $C$ is an $S$-code then its distance distribution certainly satisfied the constraints $(1), (2)$. So if the maximal value of the sum $\alpha_{\tau_1} + \ldots + \alpha_{\tau_s}$ that can be attained is $A^*$, then $A({\bf M},S)\le 1+A^*.$ (The extra 1 arises because the term $\alpha_{\tau_0}=1$ doesn't occur in this sum.)

\medskip

\noindent{\em Dual problem (LPD):} Choose a natural number $N$ and real numbers  $f_1, \ldots, f_N$
 so as to
$$
\mbox{minimize }\; f_{1} + \ldots + f_N
$$
subject to
$$
f_k\ge 0, \; k=1, \ldots, N, \quad \sum\limits_{k=1}^N {f_k \, \Phi_k(t)}\le -1, \; \,  t\in S .
$$

Thus, we have
\begin{theorem} If $A^*$ is the optimal solution to either of the primal or dual problems, then
$A({\bf M},S)\le 1+A^*.$
\end{theorem}

\medskip

\noindent {\bf 5-D. Sylvester's theorem.} Sylvester's theorem (see \cite{Pras,Mus4}) gives an answer for the following question: Suppose we know for complex numbers $t_1, \ldots, t_n$ only its power sums $ s_k$, where  $ s_k$ are real numbers. How to determine the number of real  $t_i$ in an interval $[a,b]$?

Let
$$
R_m:=
\left(
\begin{array}{cccc}
 s_0 &  s_1
&\ldots &  s_{m-1}\\
 s_1 &  s_2
&\ldots &  s_m\\
\vdots & \vdots
&\ddots & \vdots\\
 s_{m-1} &  s_m
&\ldots &  s_{2m-2}
\end{array}
\right),
$$
$$
F_m^+(a):=
\left(
\begin{array}{cccc}
s_1-as_0 &  s_2-as_1
&\ldots &  s_m-as_{m-1}\\
s_2-as_1 & s_3-as_2
&\ldots &  s_{m+1}-as_m\\
\vdots & \vdots
&\ddots & \vdots\\
s_{m}-as_{m-1} &  s_{m+1}-as_m
&\ldots &  s_{2m-1}-as_{2m-2}
\end{array}
\right),
$$
$$
F_m^-(b):=
\left(
\begin{array}{cccc}
bs_0-s_1 &  bs_1-s_2
&\ldots &  bs_{m-1}-s_m\\
bs_1-s_2 & bs_2-s_3
&\ldots &  bs_{m}-s_{m+1}\\
\vdots & \vdots
&\ddots & \vdots\\
bs_{m-1}-s_{m}&  bs_{m}-s_{m+1}
&\ldots &  bs_{2m-2}-s_{2m-1}
\end{array}
\right),
$$
and
$$
H_m(a,b)=H_m(s_0,s_1\ldots,s_{2m-1},[a,b]):=
\left(
\begin{array}{cccc}
R_m & 0 & 0\\
0& F_m^+(a) &  0\\

0 & 0 & F_m^-(b)
\end{array}
\right).
$$

\begin{theorem} Consider real numbers $t_1,\ldots,t_n$ in $[a,b]$. Then for any natural number $m$ the matrix $H_m(a,b)$ is positive semidefinite.
Moreover, for  $m\ge\rk(R_n)$ the converse holds, i.e. if $H_m(a,b)\succeq 0$, then  complex numbers  $t_1,\ldots,t_n$ with real $s_k$ are real numbers in $[a,b]$.
\end{theorem}

In fact, the constraint $H_m(a,b)\succeq 0$ doesn't depend on $n=s_0$. Indeed, let
$$
\bar s_k:=\frac{t_1^k+\ldots+t_n^k}{n}=\frac{s_k}{n}, \quad \bar H_m(a,b)=\frac{1}{n} H_m(a,b).
$$
In other words, $ \bar H_m(a,b)$ can be obtained by substituting $\bar s_k$ for $s_k$ in
$H_m(a,b)$:
$$
\bar H_m(a,b)= \bar H_m(\bar s_1,\ldots, \bar s_{2m-1},[a,b]):=H_m(1,\bar s_1,\ldots,\bar s_{2m-1},[a,b]).
$$

Thus $H_m(a,b)\succeq 0$ if and only if $\bar H_m(a,b)\succeq 0$.

\medskip

\noindent {\bf 5-E. Semidefinite programming.}
The standard form of the SDP problem is the following \cite{Todd,VB}:

\medskip

\noindent {\em Primal Problem}: $$ \mbox{minimize } \; \, c_1x_1+\ldots+c_{\ell}x_\ell$$
subject to
$$X\succeq 0, \; \mbox{ where }\; X=T_1x_1+\ldots+T_{\ell}x_{\ell}-T_0.$$

\noindent {\em Dual Problem}: $$ \mbox{maximize } \; \, \langle T_0,Y\rangle$$
subject to
$$\langle T_i,Y\rangle=c_i,\; i=1,\ldots {\ell},$$ $$ Y\succeq 0.$$

Here $T_0, T_1, \ldots, T_{\ell}, X,$ and $Y$ are real  $N\times N$ symmetric matrices,
$(c_1,\ldots,c_{\ell})$ is a cost vector, $(x_1, \ldots, x_{\ell})$ is a variable vector, and by
$\langle A,B\rangle$ we denote the inner product, i.e. $\langle A,B\rangle=\tr(AB)=\sum{a_{ij}b_{ij}}$.

\medskip

\noindent {\bf 5-F. The SDP bounds.} From here on we assume that $\Phi_k(t)$ is a polynomial of degree $k, \;
\Phi_k(\tau_0)=1$, and $S=T\cap [a,b]$ (the most interesting case for coding theory and sphere packings).

In fact, the optimal solution $A^*$ of the ${LPP}$ and  ${LPD}$ problems in 5-C depends only on the family of polynomials $\Phi:=\{\Phi_k(t)\}_{k=0,1,\ldots}$ and $[a,b]$. We denote $1+A^*$ by
$LP(\Phi,[a,b])$.

Since
$$\Phi_k(t)=p_{k0}+p_{k1}t+\ldots+p_{kk}t^k,$$
we have
$$
\Phi_k(t_1)+\ldots+\Phi_k(t_\ell)=\sum\limits_{d=0}^k{p_{kd}s_d}, \quad s_d=t_1^d+\ldots+t_\ell^d.
$$

Let $C=\{v_i\}$ be an $S$-code on ${\bf M}$, and let $\tau_{i,j}=\tau(v_i,v_j)$. Note that the number of ordered pairs $(v_i,v_j), \; i\ne j$, equals $n=|C|(|C|-1)$.  Then $(5.1)$ can be written in the form:
$$
y+p_{k0}+\sum\limits_{d=1}^k {p_{kd}x_d} \ge 0, \eqno (5.3)$$ where
$$
y=\frac{1}{|C|-1},  \quad
x_d=\bar s_d=\frac{s_d}{n}, \quad
s_d=\sum\limits_{i,j, i\ne j} {\tau_{i,j}^d}.
$$

From  Theorem 5.4  we have for any  $m$:
$$
\bar H_m(x_1,\ldots,x_{2m-1},[a,b])\succeq 0. \eqno (5.4)
$$

Now we introduce the simplest SDP bound.

\medskip

\noindent {$SDP_0$ {\em  Problem:}   Choose a natural number $m$ and real numbers $y, x_1, \ldots, x_{2m-1}$ so as to

$$\mbox{minimize } \; \, y$$
subject to
$$
y+\sum\limits_{i=1}^k {p_{ki}x_i} \ge -p_{k0}, \quad k=1,\ldots,2m-1,
$$
$$
\bar H_m(x_1,\ldots,x_{2m-1},[a,b])\succeq 0.
$$

\medskip

Note that in $(5.3)$: $|C|=(1+y)/y.$ Thus
\begin{theorem} If $y^*$ is the optimal solution of the $SDP_0$ problem, then
$$A({\bf M},S)\le SDP_0(\Phi,[a,b],m):=\frac {1+y^*}{y^*}.$$
\end{theorem}

Since we just substituted
 $H_m\succeq 0$ for
$t \in S$ in the $LPP$ problem, we can expect that $\; SDP_0(\Phi,[a,b])=LP(\Phi,[a,b])$ (see details in \cite{Mus4}).

In fact, for a continuous ${\bf M}$ the $LPP$ and $LPD$ problems are not finite linear programming problems.  These problems can be solved only via discretization. For instance, Odlyzko and Sloane \cite{OdS} ($\equiv$ \cite[Chapter 13]{CS}) applied $LPD$ for upper bounds on kissing numbers, where they replaced $S$ by 1001 equidistant points in $S$. For the $LPP$ problem it is not clear how to do a discretization of $\{\alpha_\tau\}$.
On the other hand, for a given $m$ the $SDP_0$ is a finite primal SDP problem. As a by-product of  solutions of this problem we have bounds on $|C|$ and  power sums $s_k$
(see  \cite[Section 5]{Mus4}).

In \cite[Section 6]{Mus4} is shown that some recent extensions of Delsarte's method can be reformulated as SDP problems ($SDPA$). Section 7 extends the $SDPA$ bounds to subsets of a 2-point-homogeneous space and shows that some upper bounds for codes can be improved. In particular we obtain new bounds of one-sided kissing numbers.

\section{Multivariate positive definite functions}

\noindent {\bf 6-A. Schrijver's method.}
Recently, Schrijver \cite{Schr} using  semidefinite programming (SDP) improved some upper bounds on binary codes.
Even more recently, Schrijver's method has been adapted
for  non-binary codes (Gijswijt,  Schrijver, and Tanaka \cite{GST}), and
for spherical codes  (Bachoc and Vallentin \cite{BV,BV2,BV3}). In fact, this method using the stabilizer subgroup of the isometry group derives new positive semidefinite constraints which are stronger than linear inequalities in the Delsarte linear programming method.
We consider and extend this method for spherical codes in \cite{Mus5}. Note that this approach is different from the method considered in Section 5.


\medskip

\noindent {\bf 6-B. Multivariate positive definite functions on spheres.}
Let us consider the following problem: {\em For given points $Q=\{q_1,\ldots,q_m\}$ in  $M$ to describe the class  of  continuous functions $F(t,{\bf u},{\bf v})$ in $2m+1$ variables with $ t\in {\Bbb R},\; {\bf u},{\bf v}\in {\Bbb R}^m, \; F(t,{\bf u},{\bf v})=F(t,{\bf v},{\bf u})$ such that for arbitrary points $p_1,\ldots,p_r$ in $M$ the matrix
$
\bigl(F(t_{ij},{\bf u}_i,{\bf u}_j)\bigr)\succeq 0, \; \mbox { where } \; t_{ij}=\tau(p_i,p_j), \;
{\bf u}_i=(\tau(p_i,q_1),\ldots,\tau(p_i,q_m)).\; $}

Denote this class by $\psd(M,Q)$. If
$Q=\emptyset$,  then $\psd(M,Q)$ is the class of p.d. functions in $M$.
In this case an answer is given by the Bochner\,-\,Schoenberg theorem.

\medskip

Let $0\le m\le n-2, \; t\in {\Bbb R},\; {\bf u},{\bf v}\in{\Bbb R}^m$ for $m>0$, and
${\bf u}={\bf v}=0$ for $m=0.$
 Then the following polynomial  in $2m+1$ variables of degree $k$ in $t$ is well defined:
$$
G_k^{(n,m)}(t,{\bf u},{\bf v}):=(1-|{\bf u}|^2)^{k/2}\,(1-|{\bf v}|^2)^{k/2}\,
G_k^{(n-m)}\left(\frac{t-\langle{\bf u},{\bf v}\rangle}
{\sqrt{(1-|{\bf u}|^2)(1-|{\bf v}|^2)}}\right).
$$

In \cite{Mus5} we proved the following theorem:
\begin{theorem}
Let $0\le m\le n-2$. Let $Q=\{q_1,\ldots,q_m\}\subset{\Bbb S}^{n-1}$ with $\rk(Q)=m$.
Let $e_1, \ldots, e_m$ be an orthonormal basis
  of the linear space
 with the basis $q_1,\ldots,q_m$, and
let $L_Q$ denotes the linear transformation of coordinates.
Then
 $F\in\psd({\Bbb S}^{n-1},Q)$ if and only if
$$
F(t,{\bf u},{\bf v})=\sum\limits_{k=0}^\infty{f_k({\bf u},{\bf v})\,G_{k}^{(n,m)}(t,L_Q({\bf u}),L_Q({\bf v}))},$$
where    $f_k({\bf u},{\bf v})\succeq0$ (i.e. $f_k({\bf u},{\bf v})=h_1({\bf u})h_1({\bf v})+\ldots+h_\ell({\bf u})h_\ell({\bf v})$) for all $k\ge 0.$
\end{theorem}

Note that for the case $m=0$ that is Schoenberg's theorem \cite{Scho}, and for $m=1$ it is the Bachoc-Vallentin theorem \cite{BV2}.

\medskip

\noindent {\bf 6-C. Upper bounds for spherical codes.}
In this subsection we set up upper  bounds for spherical codes which are based on multivariate p.d. functions. These bounds extend the famous Delsarte's bound. Note that for the case  $m=1$ this bound is the Bachoc - Vallentin bound \cite{BV3}.

\begin{defn} Consider a vector $J=(j_1,\ldots,j_d)$. Split the set of numbers $\{j_1,\ldots,j_d\}$ into maximal subsets $I_1,\ldots,I_k$ with equal elements. That means,  if $I_r=\{j_{r_1},\ldots,j_{r_s}\}$, then  $j_{r_1}=\ldots=j_{r_s}=a_r$ and all other $j_\ell\ne a_r$.
Without loss of generality it can be assumed that $i_1=|I_1|\ge\ldots\ge i_k=|I_k|>0$.   (Note that we have $i_1+\ldots+i_k=d$.) Denote by $\psi(J)$ the vector $\omega=(i_1,\ldots,i_k)$.

Let
$$
W_d:=\{\omega=(i_1,\ldots,i_k): i_1+\ldots+i_k=d, \; i_1\ge\ldots\ge i_k>0, \; i_1,\ldots,i_k\in {\Bbb Z}\}.
$$
Let $\omega\in W_d$. Denote
$$
\tilde q_\omega(N):=\#\{J=(j_1,\ldots,j_d)\in\{1,\ldots,N\}^d: \psi(J)=\omega\},
$$
$$
q_\omega(N):=\frac{\tilde q_\omega(N)}{N}.
$$
\end{defn}

It is not hard to see that if  $\omega\in W_d$, then $q_\omega(N)$ is a polynomial of degree $d-1$,  and  
$$
\sum\limits_{\omega\in W_d} {q_\omega(N)} = N^{d-1}.
$$

\begin{defn}
 For any vector ${\bf x}=\{x_{ij}\}$ with $1\le i<j\le d$ denote by
$A({\bf x})$  a symmetric $d\times d$ matrix $\bigl(a_{ij}\bigr)$  with all $a_{ii}=1$ and  $a_{ji}=a_{ij}=x_{ij}, \; i<j$.

Let   $\; 0<\theta<\pi$ and
$$
X(\theta):=\{{\bf x}=\{x_{ij}\}: x_{ij}\in [-1,\cos{\theta}] \mbox{ or } x_{ij}=1, \; 1\le i<j\le d\}.
$$

Now for any ${\bf x}=\{x_{ij}\}\in X(\theta)$ we  define a vector $J({\bf x})=(j_1,\ldots,j_d)$ such that $j_k=k$ if  there are no $i<k$ with  $x_{ik}=1$, otherwise $j_k=i$, where $i$ is the minimum index with $x_{ik}=1$.

Let $\omega\in W_d$. Denote
$$
D_\omega(\theta):=\{{\bf x}\in X(\theta): \psi(J({\bf x}))=\omega  \; \mbox{ and } \; A({\bf x})\succeq0\}.
$$
Let $f({\bf x})$ be a real function in ${\bf x}$, and let
$$
{B_\omega(\theta,f)}:=\sup\limits_{{\bf x}\in D_\omega(\theta)} {f({\bf x})}.
$$
\end{defn}

Note that the assumption  $A({\bf x})\succeq0$ implies  existence of unit vectors $p_1,\ldots,p_d$ such that $A({\bf x})$ is the Gram matrix of these vectors, i.e. $x_{ij}=\langle p_i,p_j\rangle$. Moreover, if $x_{ij}=1$, then $p_i=p_j$. In particular, $D_{(d)}(\theta)=\{(1,\ldots,1)\}$ and therefore ${B_{(d)}(\theta,f)}=f(1,\ldots,1)$.

\begin{defn}
Let ${\bf x}=\{x_{ij}\}$, where $1\le i<j\le m+2\le n$, and let $A({\bf x})\succeq0$.
Then there exist $P=\{p_1,\ldots,p_{m+2}\}\subset{\Bbb S}\sp {n-1}$  such that $x_{ij}=\langle p_i,p_j\rangle$.
Let $F({\bf x})$ be a continuous function in ${\bf x}$ with $F({\bf \tilde x}_{k\ell})=F({\bf x})$ for all ${\bf \tilde x}_{k\ell}$ that 
can be obtained by interchanging two points $p_k$ and $p_\ell$ in $P$.  We say that $F({\bf x})\in\psd^n_m$ if for all ${\bf x}$ with  $A({\bf x})\succeq0$ we have  $\tilde F(x_{12},{\bf u_1},{\bf u_2})\in\psd({\Bbb S}\sp {n-1},Q({\bf x}))$, where ${\bf u_i}=(x_{i3},\ldots,x_{i,m+2})$, $Q({\bf x})=\{p_3,\ldots,p_{m+2}\}$, and $\tilde F(x_{12},{\bf u_1},{\bf u_2})=F({\bf x})$.
\end{defn}


For the classical case $m=0$  Schoenberg's theorem says that $f\in \psd^n_0$ if and only if
$$f(t)=\sum\limits_{k}{f_kG_k^{(n)}(t)}$$ with all $f_k\ge 0$. It is not hard using Theorem 6.1 to describe class of functions in $\psd^n_m$ for all $m \le n-2$. 

\medskip

Let $C$ be an $N$-element subset of the unit sphere ${\Bbb S}\sp {n-1}\subset {\Bbb R}\sp n$.  It is called an    $(n,N,\theta)$ {\it spherical code} if every pair  of distinct points  $(c,c')$ of $C$ have inner product $\langle c,c'\rangle$ at most $\cos{\theta}$.


In \cite{Mus5} we proved the following theorem.
\begin{theorem}
Let $f_0>0, \; 0\le m\le n-2$, and $F({\bf x})=f({\bf x})-f_0\in\psd^n_m$. 
Then an $(n,N,\theta)$ spherical code satisfies
$$
f_0N^{m+1}\le \sum\limits_{\omega\in W_{m+2}} {B_\omega(\theta,f)\,q_\omega(N)}.
$$
\end{theorem}

It is easy to see for $m=0$ that $q_{(2)}(N)=1, \; q_{(1,1)}(N)=N-1,$ and  $B_{(2)}(\theta,f)=f(1)$. Therefore, from Theorem 6.2  we have
$$f_0N\le f(1)+B_{(1,1)}(\theta,f)(N-1).$$ Suppose $B_{(1,1)}(\theta,f)\le0$, i.e. $f(t)\le 0$ for all $t\in [-1,\cos{\theta}]$. Thus for $(n,N,\theta)$ spherical code we obtain Delsarte's bound
$$
N\le \frac{f(1)}{f_0}.
$$

\medskip

\medskip

The Bachoc-Vallentin bound \cite[Theorem 4.1]{BV3}
is the bound in Theorem 6.2 for $m=1$ and $B_{(1,1,1)}(\theta,f)\le0$. Indeed, let $B_{(2,1)}(\theta,f)\le B$. Since $q_{(3)}(N)=1, \; q_{(2,1)}(N)=3(N-1),$  and $B_{(3)}(\theta,f)=f(1,1,1)$, we have
$$
f_0N^2\le f(1,1,1)+3(N-1)B.
$$

Let us consider Theorem 6.2 also for the case $m=2$ with $B_{(1,1,1,1)}(\theta,f)\le0$. Let $B_{(3,1)}(\theta,f)\le B_1, \; B_{(2,2)}(\theta,f)\le B_2,$ and $B_{(2,1,1)}(\theta,f)\le B_3$. Then
$$
f_0N^3\le f(1,1,1,1,1,1)+4(N-1)B_1+3(N-1)B_2+6(N-1)(N-2)B_3.
$$

Let $f({\bf x})$ be a polynomial of degree  $d$.
Then the assumptions in Theorem 6.2
 can be written as  positive semidefinite constraints for the coefficients of $F$ (see for details \cite{BV,BV2,BV3,GST,Schr}).
Actually, the bound given by Theorem 6.1 can be obtained as a solution  of an SDP (semidefinite programming) optimization problem. In \cite {BV,BV2} using numerical solutions of the SDP problem for the case $m=1$
has obtained new upper bounds for the kissing numbers and for the one-sided kissing numbers in several dimensions $n\le 10$.

However, the dimension of the corresponding SDP problem is growth so fast
whenever $d$ and $m$ are increasing
that this problem can be treated numerically only for relatively small $d$ and small $m$.
It is an interesting problem to find  (explicitly) 
suitable  polynomials $F$ for Theorem 6.2 and using it to obtain new bounds for spherical codes.



\medskip

O. R. Musin, Department of Mathematics, University of Texas at Brownsville, 80 Fort Brown, Brownsville, TX, 78520.

 {\it E-mail address:} oleg.musin@utb.edu


\begin{thebibliography}{99}

\bibitem{AB1}
V.V. Arestov and A.G. Babenko, On Delsarte scheme of estimating the contact numbers,
Proc. of the Steklov Inst. of Math. {\bf 219} (1997), 36-65.

\bibitem{BV}
C. Bachoc and F. Vallentin, New upper bounds for kissing numbers from semidefinite programming,
J. Amer. Math. Soc. {\bf 21} (2008), 909-924.

\bibitem{BV2}
C. Bachoc and F. Vallentin,  Semidefinite programming, multivariate orthogonal polynomials, and codes in spherical caps, preprint, October 2006, arXiv:math.MG/0610856.

\bibitem{BV3}
C. Bachoc and F. Vallentin, Optimality and uniqueness of the (4,10,1/6) spherical code, J. Comb. Theory Ser. A {\bf 116} (2008), 195-204. 


\bibitem{BS}
E. Bannai and N.J.A. Sloane, Uniqueness of certain spherical codes, Canadian J. Math. {\bf 33} (1981), 437-449.

\bibitem{BM}
A. Barg and O.~R. Musin, Codes in spherical caps,
Advances in Mathematics of Communication, {\bf 1} (2007), 131-149.


\bibitem{Blo1} A. Blokhuis,  A new upper bound for the cardinality of  $2$-distance set in  Euclidean space, {\it Ann. Discrete Math.},  {\bf 20}  (1984), 65-66.


\bibitem{Boc}
S. Bochner, Hilbert distances and positive definite functions, Ann. Math. {\bf 42} (1941), 647-656.

\bibitem{Car}
B.~C. Carlson, Special functions of applied mathematics, Academic Press, 1977.



\bibitem{CS}
J.~H. Conway and N.~J.~A. Sloane, Sphere Packings, Lattices, and Groups, New York, Springer-Verlag, 1999 (Third Edition).

\bibitem{Del1}
Ph. Delsarte, Bounds for unrestricted codes by linear programming, Philips Res. Rep., {\bf 27}, 1972, 272-289.

\bibitem{Del2}
Ph. Delsarte, J.M. Goethals and J.J. Seidel, Spherical codes and designs, Geom. Dedic., {\bf 6}, 1977, 363-388.


\bibitem{ES}
S. J. Einhorn and I. J. Schoenberg, On Euclidean sets having only two distances between points I, II, {\it Indag. Math.}, {\bf 28} (1966), 479-488, 489-504. (Nederl. Acad. Wetensch. Proc. Ser. A69)


\bibitem{Erd}
A. Erd\'elyi, editor, Higher Transcendental Function, McGraw-Hill, NY, 3 vols, 1953, Vol. II,
Chap. XI.



\bibitem{GST}
D.~C. Gijswijt, A. Schrijver, and H. Tanaka,
New upper bounds for
nonbinary codes based on the Terwilliger algebra and semidefinite
programming, JCTA {\bf 113}(8), 2006,  p.1719--1731.


\bibitem{Kab}
G.~A. Kabatiansky and V.~I. Levenshtein, Bounds for packings on a sphere and in space,
Problems of Information Transmission, {\bf 14}(1), 1978, 1-17.




\bibitem{LRS}
D. G. Larman, C. A. Rogers, and J. J. Seidel, On two-distance sets in Euclidean space, {\it Bull. London Math. Soc.}, {\bf 9} (1977), 261-267.


\bibitem{LeS}
P. W. H. Lemmens and J. J. Seidel, Equiangular lines, {\it J. Algebra}, {\bf 24} (1973), 494-512.

\bibitem{Lev2}
V.I. Levenshtein, On bounds for packing in $n$-dimensional Euclidean space, Sov. Math. Dokl.
{\bf 20}(2), 1979, 417-421.

\bibitem{Levm}
V.I. Levenshtein, Universal bounds for codes and design, Handbook of coding theory (V. Press and W.C. Huffman, eds.), vol. 1, Elsevier Science, Amsterdam, 1998, 499-648.


\bibitem{Lis}
P. Lison\v{e}k, New maximal two-distance sets, {\it J. Comb. Theory, Ser. A}, {\bf 77} (1997), 318-338.


\bibitem{Mus}
O.~R. Musin, The problem of the twenty-five spheres, Russian Math. Surveys, {\bf 58} (2003), 794-795.

\bibitem{MusJ}
O.~R. Musin, An extension of Delsarte's method. The kissing problem in three and four dimensions, The Proceedings of COE Workshop on Sphere Packings (Nov. 1st - Nov. 5th, 2004), Kyushu University, Japan, 2005, 1-25.

\bibitem{Mus13}
O.~R. Musin, The kissing problem in three dimensions, Discrete Comput. Geom., {\bf 35} (2006), 375-384.


\bibitem{Mus3}
O.~R. Musin, The one-sided kissing number in four dimensions,
Periodica Math. Hungar., {\bf 53} (2006), 209-225.

\bibitem{Mus2}
O.~R. Musin, The kissing number in four dimensions, Ann. of Math., {\bf 168} (2008), no. 1, 1-32.

\bibitem{Mus4}
O.~R. Musin, Bounds for codes by semidefinite programming, Proc. Steklov Inst. Math. {\bf 263} (2008), 134-149.

\bibitem{Mus5}
O.~R. Musin, {Multivariable positive definite functions on spheres}, arXiv:math.MG/0701083.

 \bibitem{Mus6}
O.~R. Musin, Spherical two-distance sets, arXiv:math.MG/0801.3706, to appear in JCTA

\bibitem{Noz}
H. Nozaki, New upper bound for the cardinalities of $s$-distance sets on the unit spheres, preprint, June 2008.

\bibitem{OdS}
A.M. Odlyzko and N.J.A. Sloane, New bounds on the number of unit spheres that
can touch a unit sphere in $n$ dimensions, J. of Combinatorial Theory
A26(1979), 210-214.


\bibitem{PZ}
F. Pfender and G.~M. Ziegler, Kissing numbers, sphere packings, and some unexpected proofs, Notices Amer. Math. Soc., {\bf 51}(2004), 873-883.

\bibitem{Pras}

V.V. Prasolov, Polynomials, New York, Springer-Verlag, 2005



\bibitem{Scho}
I.~J. Schoenberg,  Positive definite functions on spheres, Duke Math. J.,
{\bf 9} (1942), 96-107.

\bibitem{Schr}
A. Schrijver, New code upper bounds from the Terwilliger algebra and semidefinite programming, IEEE Trans.  Inform. Theory {\bf 51} (2005), 2859--2866.

\bibitem{SvdW2}
K. Sch\"utte and B.L. van der Waerden, Das Problem der dreizehn Kugeln, Math. Ann. {\bf 125} (1953), 325-334.

\bibitem{Todd}
M.J. Todd, Semidefinite optimization, Acta Numerica, {\bf 10} (2001), 515-560.

\bibitem{VB}
L. Vandenberghe and S. Boyd, Semidefinite programming, SIAM Review, {\bf 38} (1996), 49-95.

\bibitem{Wang}
H.-C. Wang. Two-point homogeneous spaces, Ann. Math. {\bf 55} (1952), 177-191.

\end{thebibliography}
\end{document}